\documentclass{kms-j}

\newcommand{\publname}{}
\issueinfo{}
  {}
  {}
  {}
\pagespan{1}{}
\copyrightinfo{}
  {Korean Mathematical Society}

\usepackage{graphicx}
\allowdisplaybreaks

\theoremstyle{plain}
\newtheorem{theorem}{Theorem}[section]
\newtheorem{proposition}[theorem]{Proposition}
\newtheorem{lemma}[theorem]{Lemma}

\theoremstyle{definition}

\theoremstyle{remark}

\begin{document}

\title[Real hypersurfaces with Miao-Tam critical metrics]
{Real hypersurfaces with Miao-Tam critical metrics of complex space forms}

\author[X.M. Chen]{Xiaomin Chen}
\address{Xiaomin Chen \\ College of Science \\ China University of Petroleum-Beijing \\ Beijing 102249, China}
\email{xmchen@cup.edu.cn}
\thanks{The author is supported by the Science Foundation of China
University of Petroleum-Beijing(No.2462015YQ0604) and partially supported
by  the Personnel Training and Academic
Development Fund (2462015QZDX02).}


\begin{abstract}
Let $M$ be a real hypersurface of a complex space form with constant curvature $c$. In this paper, we study the hypersurface $M$ admitting Miao-Tam critical metric, i.e.
the induced metric $g$ on $M$ satisfies the equation:$-(\Delta_g\lambda)g+\nabla^2_g\lambda-\lambda Ric=g$, where $\lambda$ is a smooth function on $M$.
At first, for the case where $M$ is Hopf, $c=0$ and $c\neq0$ are considered respectively. For the non-Hopf case, we prove that the ruled real hypersurfaces of non-flat complex space forms do not admit Miao-Tam critical metrics.
Finally, it is proved that a compact hypersurface of a complex Euclidean space admitting Miao-Tam critical metric with $\lambda>0$ or $\lambda<0$ is a sphere and a compact hypersurface of a non-flat complex space form does not exist such a critical metric.
\end{abstract}

\maketitle

\section{Introduction}

Recall that on a compact Riemannian manifold $(M^n,g), n>2$ with a
smooth boundary $\partial M$ the metric $g$ is referred as \emph{Miao-Tam critical metric} if there
exists a smooth function $\lambda : M^n\rightarrow \mathbb{R}$ such that
\begin{equation}\label{1}
-(\Delta_g\lambda)g+\nabla^2_g\lambda-\lambda Ric = g
\end{equation}
on $M$ and $\lambda=0$ on $\partial M$, where $\Delta_g, \nabla^2_g\lambda$ are the Laplacian, Hessian operator
with respect to the metric $g$ and $Ric$ is the $(0,2)$ Ricci tensor of $g$. The
function $\lambda$ is known as the potential function. The equation \eqref{1} is called as \emph{Miao-Tam equation}.
Applying this equation,  Miao-Tam in \cite{MT} classified Einstein and conformally flat Riemannian manifolds. Particularly,  they proved that any Riemannian metric $g$ satisfying the equation
\eqref{1} must have constant scalar curvature.
 Recently, Patra-Ghosh studied the Miao-Tam equation on certain class of
odd dimensional Riemannian manifolds, namely contact metric manifolds (see
\cite{PG,PG1}). It was proved that a complete $K$-contact metric satisfying
the Miao-Tam equation is isometric to a unit sphere. Wang-Wang \cite{WW} also considered an almost Kenmotsu manifold with Miao-Tam critical metric.

An $n$-dimensional complex space form is an $n$-dimensional K\"ahler manifold with constant sectional
curvature $c$. A complete and simple connected complex space form  is complex analytically isometric to a complex projective space
 $\mathbb{C}P^n$ if $c>0$, a complex hyperbolic space $\mathbb{C}H^n$ if $c<0$, a complex Euclidean space $\mathbb{C}^n$ if $c=0.$
The complex projective and complex hyperbolic spaces are called \emph{non-flat complex space forms} and denoted by $\widetilde{M}^n(c)$.
Let $M$ be a real hypersurface of a complex space form, then
 there exists an almost contact structure $(\phi,\eta,\xi,g)$ on $M$ induced
 from the complex space form.  In particular, if $\xi$ is an eigenvector of shape operator $A$ then $M$ is
called a \emph{Hopf hypersurface}.
 Since there are no
Einstein real hypersurfaces in non-flat complex space forms(\cite{CR,M}), Cho and Kimura \cite{CK,CK2} considered a generalization of Einstein metric, called Ricci soliton, which satisfies
\begin{equation*}
 \frac{1}{2}\mathcal{L}_V g+Ric-\rho g=0,
\end{equation*}
where $V$ and $\rho$ are the potential vector field and some constant on $M$, respectively.
 They proved that a compact contact-type hypersurface with a Ricci soliton in $\mathbb{C}^n$ is a sphere and a compact Hopf hypersurface in a non-flat complex space form does not admit a Ricci soliton.

 From the Miao-Tam equation \eqref{1}, we remark that the Miao-Tam critical metric can also be viewed as a generalization of the Einstein metric since the critical metric will become an Einstein metric if the potential function $\lambda$ is constant. Thus the above results intrigue us to study the real hypersurfaces admitting  Miao-Tam critical metrics of complex space forms.
In this article, we mainly study the Hopf hypersurfaces in complex space forms as well as a class of non-Hopf hypersurfaces in non-flat complex space forms.
For a compact real hypersurface with Miao-Tam critical metric, we also get a result.

This paper is organized as follows: In Section 2 we recall some basic concepts and related results. In Section 3, we consider respectively the Hopf hypersurfaces with Miao-Tam critical metrics of non-flat complex space forms and complex Euclidean spaces, and one class of non-Hopf hypersurfaces of non-flat complex space forms
is considered in Section 4. In the last section we will prove the result of compact real hypersurfaces with Miao-Tam critical metrics.

\section{Some basic concepts and related results}

  Let ($\widetilde{M}^n,\widetilde{g})$ be a complex $n$-dimensional K\"ahler manifold
   and $M$ be an immersed, without boundary, real hypersurface of $\widetilde{M}^n$ with the induced metric $g$.
Denote by $J$ the complex structure on $\widetilde{M}^n$. There exists a local defined
unit normal vector field $N$ on $M$ and we write $\xi:=-JN$
by the structure vector field of $M$.
 An induced one-form $\eta$ is defined by
$\eta(\cdot)=\widetilde{g}(J\cdot,N)$, which is dual to $\xi$.  For any vector field $X$ on $M$ the tangent part of $JX$
is denoted by $\phi X=JX-\eta(X)N$. Moreover, the following identities hold:
\begin{equation}\label{eq2.1}
\phi^2=-Id+\eta\otimes\xi,\quad\eta\circ \phi=0,\quad\phi\circ\xi=0,\quad\eta(\xi)=1,
\end{equation}
\begin{equation}\label{eq2.2}
g(\phi X,\phi Y)=g(X,Y)-\eta(X)\eta(Y),
\end{equation}
\begin{equation}\label{eq2.3}
g(X,\xi)=\eta(X),
\end{equation}
where $X,Y\in\mathfrak{X}(M)$. By \eqref{eq2.1}-\eqref{eq2.3}, we know that $(\phi,\eta,\xi,g)$ is an almost
contact metric structure on $M$.

Denote by $\nabla, A$ the induced Riemannian connection and the shape operator on $M$, respectively.
Then the Gauss and Weigarten formulas are given by
\begin{equation}\label{eq:2.4}
\widetilde{\nabla}_XY=\nabla_XY+g(AX,Y)N,\quad\widetilde{\nabla}_XN=-AX,
\end{equation}
 where $\widetilde{\nabla}$ is the connection on $\widetilde{M}^n$ with respect to $\widetilde{g}$.
Also, we have
\begin{equation}\label{eq2.5}
  (\nabla_X\phi)Y=\eta(Y)AX-g(AX,Y)\xi,\quad\nabla_X\xi=\phi AX.
\end{equation}
In particular, $M$ is said to be a \emph{Hopf hypersurface} if the structure vector field $\xi$ is an eigenvector of $A$.

From now on we always assume that the sectional curvature of $\widetilde{M}^n$ is contant $c$. When $c=0$, $\widetilde{M}^n$ is complex Euclidean space $\mathbb{C}^n$. When $c\neq0$, $\widetilde{M}^n$ is a non-flat complex space form, denoted by $\widetilde{M}^n(c)$, then from \eqref{eq:2.4}, we know that the curvature tensor $R$ of $M$ is given by
\begin{align}\label{eq2.6}
R(X,Y)Z=\frac{c}{4}\Big(&g(Y,Z)X-g(X,Z)Y+g(\phi Y,Z)\phi X-g(\phi X,Z)\phi Y\\
&+2g(X,\phi Y)\phi Z)\Big)+g(AY,Z)AX-g(AX,Z)AY,\nonumber
\end{align}
and the shape operator $A$ satisfies
\begin{equation}\label{eq2.7}
 (\nabla_XA)Y-(\nabla_YA)X=\frac{c}{4}\Big(\eta(X)\phi Y-\eta(Y)\phi X-2g(\phi X,Y)\xi\Big)
\end{equation}
for any vector fields $X,Y,Z$ on $M$.
From \eqref{eq2.6}, we get for the Ricci tensor $Q$ of type $(1,1)$:
\begin{equation}\label{eq2.8}
 QX=\frac{c}{4}\{(2n+1)X-3\eta(X)\xi\}+hAX-A^2X,
\end{equation}
where $h$ denotes the mean curvature of $M$(i.e. $h=trace(A)$). We denote $S$ the scalar curvature of $M$, i.e. $S=trace(Q).$

If $M$ is a Hopf hypersurface of $\widetilde{M}^n(c)$, $A\xi=\alpha\xi$, where $\alpha=g(A\xi,\xi).$ Due to \cite[Theorem 2.1]{NR}, $\alpha$ is constant. Remark that when $c=0$, $\alpha$ is also constant (see the proof of \cite[Lemma 1]{CK}).
Using the equation \eqref{eq2.7}, we obtain
\begin{equation}\label{eq2.12}
 (\nabla_\xi A)X=\alpha\phi AX-A\phi AX+\frac{c}{4}\phi X
\end{equation}
for any vector field $X$. Since $\nabla_\xi A$ is self-adjoint, by taking the anti-symmetry part of \eqref{eq2.12}, we get the relation:
\begin{equation}\label{2.17}
2A\phi AX-\frac{c}{2}\phi X=\alpha(\phi A+A\phi)X.
\end{equation}

As the tangent bundle $TM$ can be decomposed as $TM=\mathbb{R}\xi\oplus\mathfrak{D}$, where $\mathfrak{D}=\{X\in TM:X\bot\xi\}$, the condition $A\xi=\alpha\xi$ implies $A\mathfrak{D}\subset\mathfrak{D}$, thus we can pick up $X\in\mathfrak{D}$ such that $AX=fX$ for some function $f$ on $M$.
Then from \eqref{2.17} we obtain
\begin{equation}\label{2.12*}
(2f-\alpha)A\phi X=\Big(f\alpha+\frac{c}{2}\Big)\phi X.
\end{equation}
If $2f=\alpha$ then $c=-4f^2,$ which show that $M$ is locally congruent a horosphere in $\mathbb{C}H^n$(see \cite{B}).

Next we recall an important lemma for a Riemannian manifold satisfying Miao-Tam equation \eqref{1}.
\begin{lemma}[\cite{GP}]\label{L1}
Let a Riemannian manifold $(M^n,g)$ satisfies the Miao-Tam
equation. Then the curvature tensor $R$ can be expressed as
\begin{equation*}
R(X,Y )\nabla\lambda =X(\lambda)QY-Y(\lambda)QX+\lambda\{(\nabla_XQ)Y-(\nabla_YQ)X\}+X(\beta)Y-Y(\beta)X
\end{equation*}
for any vector fields $X,Y$ on $M$ and $\beta=-\frac{S\lambda+1}{n-1}$.
\end{lemma}

Applying this lemma we obtain
\begin{lemma}
For a Hopf real hypersurface $M^{2n-1}$ with Miao-Tam critical metric of a complex space form, the following equation holds:
\begin{align}\label{2.19}
  \lambda\alpha \Big[X(h)-\xi(h)\eta(X)\Big]=&\mu\Big(\xi(\lambda)\eta(X)-X(\lambda)\Big)+\alpha^2\xi(\lambda)\eta(X)-\alpha AX(\lambda),
\end{align}
where $\mu=\frac{c}{4}(2n-1)+\alpha h-\alpha^2-\frac{S}{2n-2}.$
\end{lemma}
\proof Replacing $Z$ in \eqref{eq2.6} by $\nabla\lambda$, we have
 \begin{align}\label{eq2.9}
R(X,Y)\nabla\lambda=\frac{c}{4}\Big(&Y(\lambda) X-X(\lambda)Y+\phi Y(\lambda)\phi X-\phi X(\lambda)\phi Y\\
&+2g(X,\phi Y)\phi\nabla\lambda)\Big)+AY(\lambda)AX-AX(\lambda)AY.\nonumber
\end{align}
By combining with Lemma \ref{L1}, we get
\begin{align}\label{eq2.11*}
  &X(\lambda)QY-Y(\lambda)QX+\lambda\{(\nabla_XQ)Y-(\nabla_YQ)X\} \\
 = &\Big(\frac{c}{4}-\frac{S}{2n-2}\Big)\Big(Y(\lambda) X-X(\lambda)Y\Big)+\frac{c}{4}\Big(\phi Y(\lambda)\phi X
-\phi X(\lambda)\phi Y\nonumber\\
&+2g(X,\phi Y)\phi\nabla\lambda\Big)+AY(\lambda)AX-AX(\lambda)AY.\nonumber
\end{align}

Now making use of \eqref{eq2.8}, for any vector fields $X,Y$ we first compute
\begin{align*}
  (\nabla_YQ)X=&\frac{c}{4}\{-3(\nabla_Y\eta)(X)\xi-3\eta(X)\nabla_Y\xi\}+Y(h)AX+h(\nabla_YA)X\\
  &-(\nabla_YA)AX-A(\nabla_YA)X\\
  =&-\frac{3c}{4}\{g(\phi AY,X)\xi+\eta(X)\phi AY\}+Y(h)AX+h(\nabla_YA)X\\
  &-(\nabla_YA)AX-A(\nabla_YA)X.
\end{align*}
By \eqref{eq2.7}, we thus obtain
\begin{align}\label{2.12}
  &(\nabla_XQ)Y-(\nabla_Y Q)X \\
   =&-\frac{3c}{4}\{g(\phi AX+A\phi X,Y)\xi+\eta(Y)\phi AX-\eta(X)\phi AY\}\nonumber\\
   &+X(h)AY-Y(h)AX+\frac{hc}{4}\Big(\eta(X)\phi Y-\eta(Y)\phi X-2g(\phi X,Y)\xi\Big)\nonumber\\
   &-(\nabla_XA)AY+(\nabla_YA)AX-\frac{c}{4}\Big(\eta(X)A\phi Y-\eta(Y)A\phi X-2g(\phi X,Y)A\xi\Big).\nonumber
\end{align}

Therefore, taking the product of \eqref{eq2.11*} with $\xi$ and using \eqref{2.12}, we conclude that
\begin{align}\label{2.17*}
  &-\frac{3c}{4}g(\phi AX+A\phi X,Y)+\alpha X(h)\eta(Y)-\alpha Y(h)\eta(X)\\
  &-g((\nabla_XA)AY+(\nabla_YA)AX,\xi)-\frac{hc-\alpha c}{2}g(\phi X,Y)\nonumber \\
 = &\frac{\mu}{\lambda}\Big(Y(\lambda)\eta(X)-X(\lambda)\eta(Y)\Big)\nonumber\\
 &+\frac{\alpha}{\lambda}AY(\lambda)\eta(X)-\frac{\alpha}{\lambda}AX(\lambda)\eta(Y),\nonumber
\end{align}
where $\mu=\frac{c}{4}(2n-1)+\alpha h-\alpha^2-\frac{S}{2n-2}.$ Moreover, using \eqref{2.17} we compute
\begin{align*}
  &g((\nabla_XA)AY-(\nabla_YA)AX,\xi) \\
   =&g(\frac{\alpha}{2}(\phi AX-A\phi X)-\frac{c}{4}\phi X,AY) -g(\frac{\alpha}{2}(\phi AY-A\phi Y)-\frac{c}{4}\phi Y,AX).
\end{align*}
Substituting this into \eqref{2.17*} we arrive at
\begin{align}\label{eq2.17}
  &-\frac{c+2\alpha^2}{4}(\phi AX+A\phi X)+\alpha X(h)\xi-\alpha \eta(X)\nabla h\\
  &+\frac{\alpha}{2}(A^2\phi X+\phi A^2X)-\frac{2hc-\alpha c}{4}\phi X\nonumber \\
 =&\frac{\mu}{\lambda}\Big(\eta(X)\nabla\lambda-X(\lambda)\xi\Big)
 +\frac{\alpha}{\lambda}\eta(X)A\nabla\lambda-\frac{\alpha}{\lambda}AX(\lambda)\xi.\nonumber
\end{align}
Finally, taking an inner product of \eqref{eq2.17} with $\xi$ gives \eqref{2.19}.\qed

\section{Hopf real hypersurfaces of complex space forms}
First of all, we assume $c\neq0$, i.e. $M^{2n-1}$ is a Hopf real hypersurface of non-flat complex space form $\widetilde{M}^n(c).$
We first consider $\alpha=0$, i.e. $A\xi=0$, then the relation \eqref{2.19} yields
\begin{equation}\label{2.19**}
  \Big(-\frac{S}{2n-2}+\frac{c}{4}(2n-1)\Big)\Big(\xi(\lambda)\xi-\nabla\lambda\Big)=0.
\end{equation}

 If $-\frac{S}{2n-2}+\frac{c}{4}(2n-1)=0$, i.e. $S=\frac{1}{2}c(n-1)(2n-1)$.  Then from \eqref{eq2.17}  we find
\begin{equation}\label{2.19*}
\frac{c}{4}(\phi AX+A\phi X)=0,
\end{equation}
which yields $\phi AX+A\phi X=0$ for all vector field $X$. This is contradictory with \cite[Corollary 2.12]{NR}.
Thus $S\neq\frac{c}{2}(n-1)(2n-1),$ and it follows from \eqref{2.19**} that $\nabla\lambda=\xi(\lambda)\xi.$
Differentiating this along $X$ gives
\begin{equation}\label{eq2.19}
  \nabla_X\nabla\lambda=X(\xi(\lambda))\xi+\xi(\lambda)\phi AX.
\end{equation}

On the other hand, from \eqref{1} we can obtain
\begin{equation}\label{eq2.20}
 \nabla_X\nabla\lambda=(1+\Delta\lambda)X+\lambda QX.
\end{equation}
Comparing \eqref{eq2.19} and \eqref{eq2.20}, we have
\begin{equation}\label{2.21}
X(\xi(\lambda))\xi+\xi(\lambda)\phi AX=(1+\Delta\lambda)X+\lambda QX.
\end{equation}
Moreover, by \eqref{eq2.8}, putting $X=\xi$ gives
\begin{equation}\label{2.22}
  \xi(\xi(\lambda))=1+\Delta\lambda+\frac{\lambda c}{2}(n-1).
\end{equation}

Choose a local orthonormal frame $\{e_i\}_{i=1}^{2n-1}$ such that $e_{2n-1}=\xi$ and $e_{n-1+i}=\phi e_{i}$ for $i=1,\cdots,n-1$.
Using the frame to contract over $X$ in \eqref{2.21}, we also derive that
\begin{equation*}
\xi(\xi(\lambda))=(1+\Delta\lambda)(2n-1)+\lambda S.
\end{equation*}
Comparing with \eqref{2.22}, we find
\begin{equation}\label{2.23}
(2n-2)(1+\Delta\lambda)+\lambda S=\frac{\lambda c}{2}(n-1).
\end{equation}

Furthermore, by taking the trace of Miao-Tam equation \eqref{1},
we get
\begin{equation}\label{2.25}
 (2-2n)\Delta\lambda-\lambda S=2n-1,
\end{equation}
which, together with \eqref{2.23}, yields
\begin{equation}\label{2.27}
\frac{\lambda c}{2}(n-1)+1=0.
\end{equation}
This show that $\lambda$ is constant. Thus $M$ is Einstein, but as is well-known that there are no Einstein hypersurfaces in a non-flat complex space form as in introduction, hence we immediately obtain
\begin{proposition}\label{P1}
A real hypersurface with $A\xi=0$ of a non-flat complex space form does not admit Miao-Tam critical metric.
\end{proposition}
Next we consider the case where $\alpha\neq0$.
If for every $X\in\mathfrak{D}$  such that $AX=\frac{\alpha}{2}X$, as before we know that $M$ is locally congruent a horosphere in $\mathbb{C}H^n$ and $c=-\alpha^2.$ Moreover, the mean curvature $h=n\alpha$ is constant. Then from \eqref{eq2.17} we can obtain $nc=-\frac{\alpha^2}{2}$. This implies $2n=1$. It is impossible.

Now choose $X\in\mathfrak{D}$ such that $AX=fX$ with $f\neq\frac{\alpha}{2}$, so from \eqref{eq2.17} we have
\begin{align*}
  &-\frac{c+2\alpha^2}{4}(f\phi X+\widetilde{f}\phi X)+\alpha X(h)\xi+\frac{\alpha}{2}(\widetilde{f}^2\phi X+f^2\phi X)-\frac{2hc-\alpha c}{4}\phi X\nonumber \\
 =&-\frac{\mu}{\lambda}X(\lambda)\xi-\frac{\alpha}{\lambda}AX(\lambda)\xi.
\end{align*}
Here we have used $A\phi X=\widetilde{f}\phi X$ with $\widetilde{f}=\frac{f\alpha+\frac{c}{2}}{2f-\alpha}$ followed from \eqref{2.12*}.
Since $\phi X\in\mathfrak{D}$, we further derive
\begin{equation}\label{2.18}
  -(c+2\alpha^2)(f+\widetilde{f})+2\alpha(\widetilde{f}^2+f^2)-(2hc-\alpha c)=0.
\end{equation}
Moreover, inserting $\widetilde{f}=\frac{f\alpha+\frac{c}{2}}{2f-\alpha}$ into the equation \eqref{2.18}, we have
\begin{align}\label{eq2.19*}
  &8\alpha f^4-4(c+4\alpha^2)f^3+(6\alpha c+8\alpha^3-8hc)f^2\\
  &+(8hc\alpha-4\alpha^2c-c^2)f+\alpha c^2+2\alpha^3c-2hc\alpha^2=0.\nonumber
\end{align}

Now we denote the roots of the polynomial by $f_1,f_2,f_3,f_4$, then from the relation between the roots and coefficients we obtain
\begin{align}\label{2.20}
\left\{
  \begin{array}{ll}
    &f_1+f_2+f_3+f_4=\frac{c+4\alpha^2}{2\alpha },  \\
    &f_1f_2+f_1f_3+f_1f_4+f_2f_3+f_2f_4+f_3f_4 =\frac{3\alpha c+4\alpha^3-4hc}{4\alpha},  \\
    &f_1f_2f_3+f_1f_2f_4+f_2f_3f_4 =-\frac{8hc\alpha-4\alpha^2c-c^2}{8\alpha},  \\
    &f_1f_2f_3f_4 =\frac{c^2+2\alpha^2c-2hc\alpha}{8}.
  \end{array}
\right.
\end{align}
As the proof of \cite[Lemma 4.2]{CK}, we can also get the following.
\begin{lemma}\label{L2}
 The mean curvature $h$ is constant.
\end{lemma}
Hence from \eqref{2.19} we conclude
\begin{equation*}
 A\nabla\lambda=\frac{\mu}{\alpha}\phi^2\nabla\lambda+\alpha\xi(\lambda)\xi.
\end{equation*}
By taking the inner product with the principal vector $X\in\mathfrak{D}$, we obtain
\begin{equation*}
(f+\frac{\mu}{\alpha})X(\lambda)=0.
\end{equation*}
If $X(\lambda)=0$ for all $X\in\mathfrak{D}$, then $\nabla\lambda=\xi(\lambda)\xi.$ As the proof of Proposition \ref{P1}, we see that $M$ is Einstein, which is impossible.

If $X(\lambda)\neq0$ for all $X\in\mathfrak{D}$, then $f+\frac{\mu}{\alpha}=0$, i.e. $M$ has only two distinct constant principal curvatures $\alpha,-\frac{\mu}{\alpha}$. Further, we see from \eqref{2.12*} that
\begin{equation}\label{3.30}
2f^2-2\alpha f-\frac{c}{2}=0.
\end{equation}
 Since the hypersurface $M$ has two distinct constant principle curvatures: $\alpha$ of multiplicity $1$ and $f$ of multiplicity $2n-2$,
it is easy to get that the mean curvature $h=\alpha+(2n-2)f$ and the scalar curvature $S=c(n^2-1)+2\alpha(2n-2)f+(2n-2)(2n-3)f^2.$ Thus
\begin{equation*}
  \mu=-\frac{3c}{4}+(2n-4)\alpha f-(2n-3)f^2.
\end{equation*}
Inserting this into the relation $f+\frac{\mu}{\alpha}=0$, we obtain
\begin{equation}\label{3.31}
  (2n-3)(\alpha f-f^2)=\frac{3c}{4}.
\end{equation}
Combining \eqref{3.30} with \eqref{3.31}, we find $nc=0,$ which is a contradiction.

If $X(\lambda)\neq0$ for some principle vector $X\in\mathfrak{D}$, and without loss general, we suppose $e_1(\lambda)\neq0$, then $Ae_1=-\frac{\mu}{\alpha}e_1$ and $A\phi e_1=\frac{\alpha\mu-\frac{c}{2}\alpha}{2\mu+\alpha^2}\phi e_1$.

Notice that if the hypersurface $M$ of $\mathbb{C}H^n$ has constant principal curvatures, the classification is as follows:
\begin{theorem}[\cite{B}]\label{B}
Let $M$ be a Hopf real hypersurface in $\mathbb{C}H^n(n\geq2)$ with constant principal curvatures. Then
$M$ is locally congruent to the following:
\begin{enumerate}
 \item $A_2$: Tubes around a totally geodesic $\mathbb{C}H^{n-1}\subset\mathbb{C}H^n$.
  \item $B$: Tubes of radius $r$ around a totally geodesic real hyperbolic space $\mathbb{R}H^n\subset\mathbb{C}H^n$.
\item $N$: Horospheres in $\mathbb{C}H^n$.
 \end{enumerate}
\end{theorem}
Since the horospheres have two distinct principal curvatures, it is impossible.
By Theorem 3.9 and 3.12 in \cite{NR}, the Type $A_2,B$ hypersurfaces have three distinct principal curvatures: $\lambda_1=\frac{1}{r}\tanh(u),\lambda_2=\frac{1}{r}\coth(u)$ and $\alpha=\frac{2}{r}\tanh(2u).$ Then $h=\alpha+(n-1)(\lambda_1+\lambda_2)=\alpha+\frac{2(n-1)}{r}\coth(2u).$ On the other hand, from Corollary 2.3(ii) in \cite{NR}, we also have
$\frac{1}{r^2}=\frac{\lambda_1+\lambda_2}{2}\alpha+\frac{c}{4}$, i.e. $c=-\frac{4}{r^2}.$
This implies from the last relation in \eqref{2.20} that
\begin{align*}
&\frac{1}{r^4}=\frac{c^2+2\alpha^2c-2hc\alpha}{8}=\frac{4n-2}{r^4}.
\end{align*}
Thus $n=\frac{3}{4}$, that is impossible.

For the case of $\mathbb{C}P^n$, the classification is as follow:
\begin{theorem}[\cite{K2,T1}]
Let $M$ be a Hopf hypersurface in $\mathbb{C}P^n(n\geq2)$ with constant principal curvatures.
 Then $M$ is an open part of
\begin{enumerate}
  \item $A_2$: a tuber over a totally geodesic complex projective space $\mathbb{C}P^k$ of radius $\frac{\pi r}{4}$ for $0\leq k\leq n-1$, where $r=\frac{2}{\sqrt{c}}$, or
  \item B: a tuber over a complex quadric $Q_{n-1}$ and $\mathbb{R}P^n$, or
  \item C: a tube around the Segre embedding of $\mathbb{C}P^1\times \mathbb{C}P^k$ into $\mathbb{C}P^{2k+1}$ for some $k\leq2$, or
\item D: a tube around the Pl\"ucker embedding into $\mathbb{C}P^9$ of the complex Grassmann manifold
$G_2(\mathbb{C}^5)$ of complex $2$-planes in $\mathbb{C}^5$, or
\item E: a tube around the half spin embedding into $\mathbb{C}P^{15}$ of the Hermitian symmetric space
$SO(10)=U(5)$.
\end{enumerate}
\end{theorem}
The Type $A_2$ and $B$ hypersurfaces have three distinct principal curvatures: $\lambda_1=-\frac{1}{r}\cot(u),\lambda_2=\frac{1}{r}\tan(u),\alpha=\frac{2}{r}\tan(2u)$(see \cite[Theorem 3.14, 3.15]{NR}). From the  first relation of \eqref{2.20}, we have
\begin{align*}
    \lambda_1+\lambda_2=\frac{c+4\alpha^2}{4\alpha }\Rightarrow-\frac{16}{r^2}=c+4\alpha^2.
\end{align*}
It gives a contradiction since $c>0$.

For the Type $C,D$ and $E$ hypersurfaces, they have five distinct principal curvatures(see \cite[Theorem 3.16, 3.17, 3.18]{NR}).
We compute
\begin{equation*}
  \frac{1}{r}\Big(-\cot(u)+\tan(u)+\cot(\frac{\pi}{4}-u)+\cot(\frac{3\pi}{4}-u)\Big)=\frac{2}{r}(1+\cot^2(2u)).
\end{equation*}
Thus the  first relation of \eqref{2.20} implies
\begin{equation*}
-\frac{24}{r^2}\cot^2(2u)=c+\frac{8}{r^2}.
\end{equation*}
It is impossible since $c>0$. So the hypersurfaces of type $C,D,E$ do not admit Miao-Tam critical metrics.

Summarizing the above discussion,  we thus assert the following:
\begin{proposition}\label{P2}
A real hypersurface with $A\xi=\alpha\xi,\alpha\neq0$ in a non-flat complex space form does not admit Miao-Tam critical metric.
\end{proposition}
Together Proposition \ref{P1} with Proposition \ref{P2}, we prove
\begin{theorem}
There exist no Hopf real hypersurfaces with Miao-Tam critical metric in non-flat complex space forms.
\end{theorem}

In the following we always assume $c=0$. That is to say that $M$ is a real hypersurface of complex Euclidean space $\mathbb{C}^n.$ First of all, if $A\xi=0$, we obtain from \eqref{2.19**}
\begin{equation*}
 S\Big(\xi(\lambda)\xi-\nabla\lambda\Big)=0.
\end{equation*}
If $S\neq0$, we have $\nabla\lambda=\xi(\lambda)\xi$. As before we can also lead to \eqref{2.27}, but it yields a contradiction since $c=0$.
Thus the scalar curvature $S=0$, and the relation \eqref{2.25} implies $\Delta\lambda=-\frac{2n-1}{2n-2}$. Actually, $\lambda=-\frac{2n-1}{4n-4}|x|^2$ on $\mathbb{R}^{2n-1}.$ Since $R(\xi,X,\xi,X)=0$ for all $X$, the sectional curvature of $M$ is also zero. By Hartman and Nirenberg's theorem in \cite{HN},
$M$ is a hyperplane or a cylinder, hence we have the following:
\begin{theorem}
Let $M^{2n-1}$ be a real hypersurface with $A\xi=0$ of complex Euclidean space $\mathbb{C}^n$. If $M$ admits Miao-Tam critical metric, it is a generalized cylinder $\mathbb{R}^{2n-1-p}\times \mathbb{S}^{p}$ or $\mathbb{R}^{2n-1}$.
\end{theorem}

When $\alpha\neq0.$ Let us choose $X\in\mathfrak{D}$ such that $AX=\beta X$ for a smooth function $\beta$, then we know $\beta\neq\frac{\alpha}{2}$, otherwise, if $\beta=\frac{\alpha}{2}$, then
$-4\beta^2=c=0$ from \eqref{2.12*}, i.e. $\beta=0.$ This is a contradiction with $\alpha\neq0$. Further, from \eqref{2.12*} we have
\begin{equation}\label{2.32}
  A\phi X=\frac{\beta\alpha}{2\beta-\alpha}\phi X.
\end{equation}
Therefore we find that the equation \eqref{eq2.19*} holds, and for $c=0$ and $f=\beta$ it becomes
\begin{align*}
  &(\beta^2-\alpha\beta)^2=0.
\end{align*}
So $\beta^2=\alpha\beta$, that means that $\beta$ is constant and further $h$ is also constant. If $\alpha=\beta$, from \eqref{2.32} we see that the shape operator can be expressed as $A=\alpha I$, where $I$ denotes the identity map. In this case, $M$ is locally congruent to a sphere.

If $\beta=0$, $A=\alpha\eta\otimes\xi$, as the proof of \cite[Theorem 1.1]{LXZ}, we know that $M$ is $\mathbb{S}^1\times\mathbb{R}^{2n-2}$. Therefore we assert the following:

\begin{theorem}
Let $M^{2n-1}$ be a real hypersurface with $A\xi=\alpha\xi,\alpha\neq0,$ of complex Euclidean space $\mathbb{C}^n$. If $M$ admits Miao-Tam critical metric, it is locally congruent to a sphere, or $\mathbb{S}^1\times\mathbb{R}^{2n-2}$.
\end{theorem}

\section{Ruled hypersurfaces of non-flat complex space forms}
In this section we study a class of non-Hopf hypersurfaces with Miao-Tam critical metric of non-flat complex space forms.
Let $\gamma:I\rightarrow\widetilde{M}^n(c)$ be any regular curve. For $t\in I$, let $\widetilde{M}^n_{(t)}(c)$ be a totally geodesic complex hypersurface
through the point $\gamma(t)$ which is orthogonal to the holomorphic plane spanned by $\gamma'(t)$ and $J\gamma'(t)$.
Write $M=\{\widetilde{M}^n_{(t)}(c):t\in I\}$. Such a construction asserts that $M$ is a real hypersurface of $\widetilde{M}^n(c)$, which is called a \emph{ruled hypersurface}. It is well-known that the shape operator $A$ of $M$ is written as:
\begin{align*}
A\xi=&\alpha\xi+\beta W(\beta\neq0), \\
AW = &\beta\xi,\\
AZ = &0\; \text{for any}\; Z\bot\xi,W,
\end{align*}
where $W$ is a unit vector field orthogonal to $\xi$, and $\alpha,\beta$ are differentiable functions on $M$. From \eqref{eq2.8}, we have
\begin{align}
Q\xi =& (\frac{1}{2}(n-1)c-\beta^2 )\xi,\label{4.32}\\
QW =& (\frac{1}{4}(2n+1)c-\beta^2 )W,\label{4.33*}\\
QZ =& (\frac{1}{4}(2n+1)c)Z\; \text{for any}\; Z\bot\xi,W.\label{4.34}
\end{align}
From these equations we know the scalar curvature $S=(n^2-1)c-2\beta^2$. Since $S$ is constant, this shows that $\beta$ is also constant. Further, the following relation
$\nabla\beta=(\beta^2+c/4)\phi W$ is valid (see \cite{K}), which yields
\begin{equation}\label{3.32}
  \beta^2+c/4=0\quad\text{and}\quad S=-(4n^2-2)\beta^2.
\end{equation}
Further, the following lemma holds:
\begin{lemma}[\cite{K}]\label{L5}
For all $Z\in\{X\in TM:\eta(X)=g(X,W)=g(X,\phi W)=0\}$, we have the following relations:
\begin{align*}
\nabla_W\phi W &= -2\beta W, \quad\nabla_WW =(\beta+\beta^2)\phi W,\\
\nabla_Z\phi W &= -\beta Z, \quad\nabla_ZW =\beta\phi Z,\\
\nabla_{\phi W}& \phi W = 0.
\end{align*}
\end{lemma}

Now putting $Y=\xi$ and $X=W$ in \eqref{eq2.11*} yields
\begin{align}\label{4.33}
  &W(\lambda)(\frac{1}{2}(n-1)c-\beta^2 )\xi-\xi(\lambda)(\frac{1}{4}(2n + 1)c-\beta^2 )W\\
  &+\lambda\{(\nabla_WQ)\xi-(\nabla_\xi Q)W\}\nonumber \\
 = &\Big(\frac{c}{4}-\frac{S}{2n-2}\Big)\Big(\xi(\lambda) W-W(\lambda)\xi\Big)+A\xi(\lambda)AW-AW(\lambda)A\xi.\nonumber
\end{align}
Because $\beta$ is constant,  from \eqref{4.33*} and \eqref{4.32}, by Lemma \ref{L5} we compute
\begin{align*}
  (\nabla_WQ)\xi-(\nabla_\xi Q)W & =\nabla_W(Q\xi)-Q\nabla_W\xi-\nabla_\xi(QW)+Q\nabla_\xi W \\
   & =-W(\beta^2)\xi+\xi(\beta^2)W=0.
\end{align*}
Inserting this into \eqref{4.33}, we conclude that
\begin{equation}\label{3.36}
\left\{
  \begin{array}{ll}
    W(\lambda)\Big[(\frac{1}{4}(2n-1)c-2\beta^2-\frac{S}{2n-2}\Big] &=0,  \\
    \xi(\lambda)\Big[\frac{1}{2}(n+1)c-2\beta^2-\frac{S}{2n-2}\Big]&= 0.
  \end{array}
\right.
\end{equation}
From \eqref{3.36}, we get $\xi(\lambda)=W(\lambda)=0$ since $\frac{1}{2}(n+1)c-2\beta^2-\frac{S}{2n-2}\neq0$, which is followed from \eqref{3.32}.

Putting $Y=\xi$ and $X=Z$ in \eqref{eq2.11*}, we have
\begin{align}\label{40}
  &Z(\lambda)(\frac{1}{2}(n-1)c-\beta^2 )\xi-\xi(\lambda)(\frac{1}{4}(2n+1)c)Z+\lambda\{(\nabla_ZQ)\xi-(\nabla_\xi Q)Z\} \\
 = &\Big(\frac{c}{4}-\frac{S}{2n-2}\Big)\Big(\xi(\lambda) Z-Z(\lambda)\xi\Big).\nonumber
\end{align}
By Lemma \ref{L5}, we also obtain
\begin{equation*}
(\nabla_ZQ)\xi-(\nabla_\xi Q)Z=-Z(\beta^2)\xi+\xi(\beta^2)Z=0.
\end{equation*}

Since $\xi(\lambda)=0$, the relation \eqref{40} becomes
\begin{equation*}
 Z(\lambda)\Big[\frac{1}{4}(2n-1)c-\beta^2-\frac{S}{2n-2}\Big]=0.
\end{equation*}
Thus $Z(\lambda)=0$ since $\frac{1}{4}(2n-1)c-\beta^2-\frac{S}{2n-2}\neq0$ as before.

By taking $X=\phi W$ and $Y=\xi$ in \eqref{eq2.11*}, a similar computation gives
\begin{align}\label{3.37}
  -\lambda\beta(\frac{1}{2}(n+2)c+\beta^2)=\Big(-\frac{S}{2n-2}+\frac{1}{4}(2n-1)c-\beta^2\Big)\phi W(\lambda).
\end{align}
Inserting \eqref{3.32} into \eqref{3.37}, we find
\begin{equation*}
\phi W(\lambda)=\frac{\lambda\beta(2n+3)(n-1)}{2n-1}.
\end{equation*}
Consequently, we obtain
\begin{equation}\label{4.37}
  \nabla\lambda=\frac{\lambda\beta(2n+3)(n-1)}{2n-1}\phi W.
\end{equation}

On the other hand, as we known $\nabla_X\nabla\lambda=\lambda QX+(1+\Delta\lambda)X$ by Miao-Tam equation \eqref{1}. When $X=Z$ and $W$ respectively, by Lemma \ref{L5} it follows respectively
from \eqref{4.33*}, \eqref{4.34} and \eqref{4.37} that
\begin{align*}
  -\frac{\lambda\beta^2(2n+3)(n-1)}{2n-1}=-\lambda(2n+1)\beta^2+(1+\Delta\lambda),\\
-2\frac{\lambda\beta^2(2n+3)(n-1)}{2n-1}=-\lambda(2n+2)\beta^2+(1+\Delta\lambda).
\end{align*}
It will give $\lambda\beta^2=0$, which is a contradiction with $\lambda,\beta\neq0$. Hence the following theorem is proved.
\begin{theorem}
There exist no ruled hypersurfaces with Miao-Tam critical metrics of non-flat complex space forms.
\end{theorem}

\section{Compact hypersurfaces of complex space forms}
For the case where $M$ is compact, we immediately obtain the following result:
\begin{theorem}
Let $M^{2n-1}$ be a compact real hypersurface admitting Miao-Tam critical metric with $\lambda>0$ or $\lambda<0$ of complex Euclidean space $\mathbb{C}^n$, then
 $M$ is a sphere. In the compact real hypersurfaces of a non-flat complex space form $\widetilde{M}^n(c)$ there does not exist such a critical metric.
\end{theorem}
\proof Write $\mathring{Ric}=Ric-\frac{S}{2n-1}g$. It is proved the following relation(see the proof of \cite[Lemma 5]{BDRR}):
\begin{equation*}
  {\rm div}(\mathring{Ric}(\nabla\lambda))=\lambda|\mathring{Ric}|^2.
\end{equation*}
Thus integrating it over $M$ gives $\mathring{Ric}=0$ if $\lambda>0$ or $\lambda<0$, that means that $Ric=\frac{S}{2n-1}g$.
Namely $M$ is Einstein. For the case of complex Euclidean space $\mathbb{C}^n$, it is proved that $M$ is a sphere,
a hyperplane, or a hypercylinder over a complete plane curve (cf. \cite{F}). But the latter two cases are not compact.
 For $c\neq0$, it is impossible since there are no Einstein hypersurfaces in a non-flat complex space form. Therefore we complete the proof.\qed

\section*{Acknowledgement}
The author would like to thank the referees for the helpful
suggestions.

\end{document}